\documentclass{lms}
\usepackage{url,graphicx,amssymb,tikz}

\newtheorem{theorem}{Theorem}[section] % 1st argument is your name for it
\newcommand{\gen}[1]{\left< #1 \right>}

\newnumbered{assertion}{Assertion}    % 1st argument is your name for it
\newnumbered{conjecture}{Conjecture}  % 2nd argument is what is printed
\newnumbered{definition}{Definition}
\newnumbered{hypothesis}{Hypothesis}
\newnumbered{remark}{Remark}
\newnumbered{note}{Note}
\newnumbered{observation}{Observation}
\newnumbered{problem}{Problem}
\newnumbered{question}{Question}
\newnumbered{algorithm}{Algorithm}
\newnumbered{example}{Example}
\newunnumbered{notation}{Notation} % This is usually unnumbered
\title[John McKay Obit]% end with percent
 {John Keith Stuart McKay: 1939 - 2022} % This is the full title of the paper
\author{Yang-Hui He}

\dedication{Non omnis moriar multaque pars mei
vitabit Libitinam - Carmina Horatii, III.30}

\classno{01Axx, 20-XX, 11-XX, 17BXX, 06Bxx, 14E16}
\extraline{Partly supported by UK STFC grant ST/J00037X/2 and Leverhulme Project Grant 2022.
}

\begin{document}
\maketitle

%\begin{abstract}
%A homage to the life and mathematics of John K.~S.~McKay.
%Obituary for the Bulletin of the London Mathematical Society. 
%\end{abstract}

%\part{Use this type of header for very long papers only}
% use lowercase except for proper names

\begin{quotation}
``I could never have approached this subject but for the untiring interest of John McKay, who has often seemed to me an emissary from some advanced galactic civilization, sent here to speed up our evolution.'' \\
\flushright{ - Jack Morava \cite{morava}}
\end{quotation}

These words of the topologist Morava summarise well the awe and, to a degree, the bemusement, which the mathematics community has felt towards John McKay. Certainly, of all mathematicians in the recent, if not entire, history of the subject, McKay was rare: he observed the vista of mathematics from afar and saw profound hidden connections, some of which were even mocked at first and seen as crazy. 

I had the honour and pleasure to be John's last close collaborator and -- as his wife Trinh quotes him -- ``his best friend'' in his old age.
For these last 10 years, I have felt I was under the wing of my childhood hero.
I never had a grandfather, both of mine long being gone before my birth. 
Yet, rather serendipitously, John McKay became a grandfather figure to me and my family, intellectually and emotionally\footnote{There were three good English friends and colleagues: John Horton Conway (1937 – 2020), Simon Phillips Norton (1952 – 2019), and John Keith Stuart McKay (1939 – 2022). Conway and Norton will feature prominently in Section 2. Alexander Masters wrote the story of Norton in {\it The Genius in my Basement} \cite{nortonbio}, Siobhan Roberts wrote the story of Conway in {\it Genius at Play} \cite{conwaybio}. I am currently writing  {\it The Genius on Skype} \cite{mckaybio}, a somewhat part three of this trilogy. 
}. When asked by the London Mathematical Society to write this obituary, I felt pride, and a tremendous sense of  duty, as well as pangs of sorrow. I only hope that I can, while recounting the ensuing story, in some small way be a Dante to his Virgil, or a Boswell to his Johnson.

\begin{figure}
\centerline{
\includegraphics[trim=0mm 0mm 0mm 0mm, clip, width=2in]{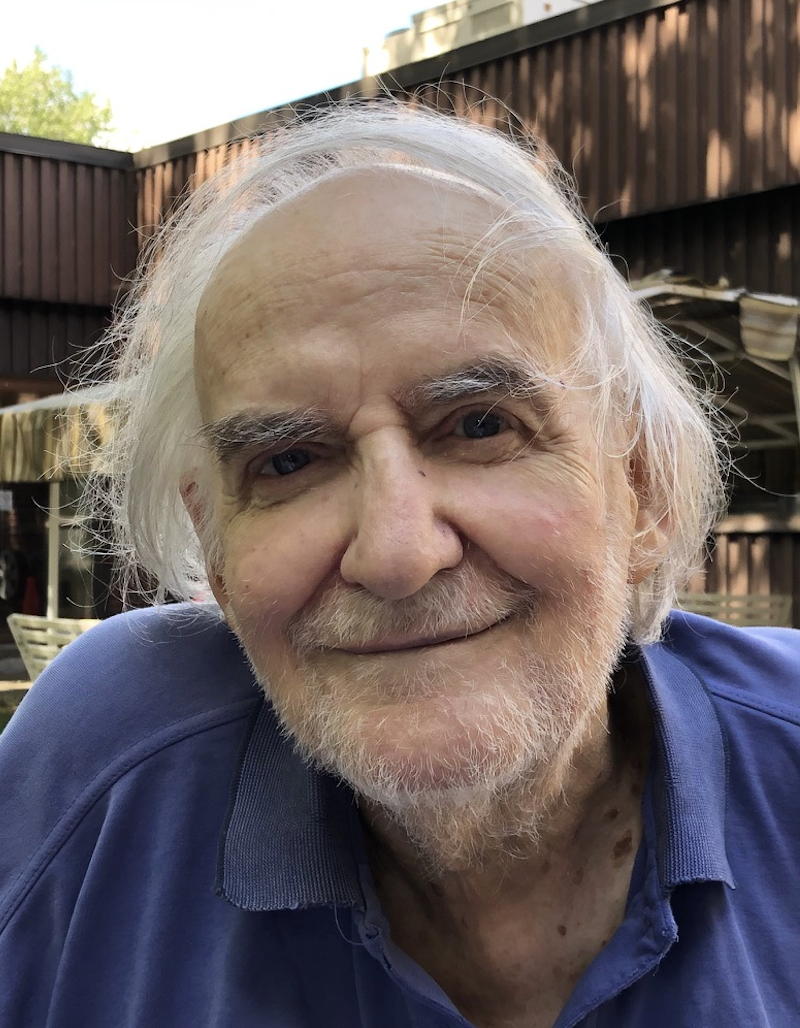}
\includegraphics[trim=0mm 0mm 0mm 0mm, clip, width=2.2in]{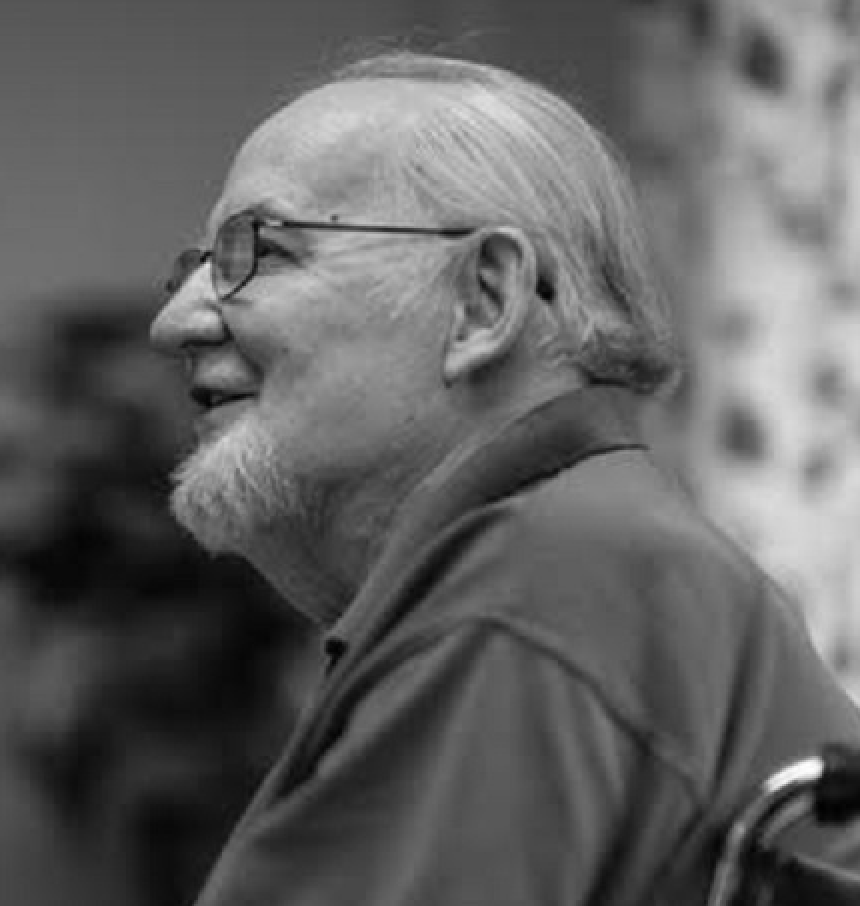}
}
\caption{John McKay, courtesy of his wife Trinh Vo-McKay}
\label{f:photoTrinh}
\end{figure}

\section{Biography} % use lowercase except for proper names
\label{bio}
John Keith Stuart McKay was born on 18 of November, 1939, in Kent, on the outskirts of London, to a comfortably-off English family. His father, Ernest Stuart McKay, ran a successful book binding business in London (the Bookbinders of London, est. 1880) and his mother, Norah Florence (n\'ee Culshaw), came from a well-connected background. To his English family and friends, McKay was known as ``Keith'' -- the ``John'' persona came later to his trans-Atlantic, and global colleagues.
He and his elder sister Elizabeth, who still lives in England, grew up in a large Victorian family home, a tranquility shaken in their early years only by the tumults of WW2, including the time when a German V-bomb landed not too far from their dwellings.
Like many boys of his generation from English families of some affluence, John was sent away to boarding school.
From 1950 to 1958, he attended the well-known public school, Dulwich College\footnote{I have been invited to speak at Dulwich a few times, both independently of and specifically about, John. Taking advantage of the friendship I have harboured with the masters, I have on many occasions suggested that a plaque be dedicated to McKay at his alma mater.}, like his father before him.

On obtaining his A-levels, McKay went to the University of Manchester to study mathematics, where he learnt about the then-nascent subject of computing.
The university was building the Manchester Atlas Supercomputer, which was then one of the most powerful machines in the world. 
This joint interest in pure mathematics, scientific computing, as well as mathematical data would have foundational implications for his future research.
Indeed, on finishing his undergraduate studies in 1962, McKay stayed on until 1964 as Nuffield Research Fellow, working with the Atlas, as one of the pioneers of the use of computers in pure mathematics.

This theme of mathematics and computers continued in 1964 when McKay was offered a PhD place at the University of Edinburgh under the joint supervision of Sidney Michaelson, Scotland's first professor of computer science, and W.~Douglas Munn, an algebraist specialising in representation theory of groups and semi-groups.
Appropriately, McKay's thesis was entitled {\em Computing with Finite Groups}, a great and first synthesis of both fields \cite{thesis}.
The Acknowledgments of the dissertation read:
``I thank my supervisors, Professor W.~D.~Munn and Professor S.~Michaelson, for much patience and many kindnesses during the preparation of this thesis. My thanks go also to my wife Wendy who has tirelessly typed and retyped the manuscript and organised me.''
From this we can deduce several facts about John.
First, he was appreciative of his mentors. Indeed, throughout his life, John was a kind and gentle man, to his teachers, peers, and the younger generation.
Second, he got married during his PhD years. This was to Wendy Sue-A-Quan, a fellow student at Edinburgh.
Third, Wendy organised him. For those who knew McKay, the opinion was ubiquitous that he was the archetype of the absent-minded professor, brilliant, soft-spoken, and charmingly disorganised.

With a doctorate -- the thesis was formally completed and submitted some years later in 1971 -- McKay became a research fellow by returning to the Atlas Computer world, this time at Chilton, just outside of Oxford, in 1967. The same year, he became a member of the London Mathematical Society.
Chilton gave McKay ample opportunity to collaborate with Oxford University, and in particular Graham Higman, then the Waynflete Chair of pure mathematics, in the explicit construction of some sporadic simple groups -- specifically, the third Janko group $J_3$ and the Held group $He$ -- with the help of the computer.
One of John's contemporaries at Chilton was John Leech, the originator of the eponymous lattice \cite{leech}, which is the unique even unimodular lattice in 24 dimensions and the best one for sphere packing therein. Leech tried in vain to get the group theory community interested in this, despite his suspicion that the enormous group of automorphisms of the lattice should be of universal interest.
Luckily, McKay got interested, then got John Horton Conway interested, and this led to Conway's
discovery and construction of three new sporadic simple groups, now called
$Co_1$, $Co_2$, and $Co_3$. The group $Co_1$, which is the quotient
of the automorphism group of the Leech lattice by its centre of size 2, 
turned out to play a crucial role in the construction of the legendary 
Monster group, with which we will spend much time shortly.

One of John's greatest powers was connecting people with ideas.
Siobhan Roberts, in her award-winning biography of Conway \cite{conwaybio}, describes it well in her playfully entertaining way: 
\begin{quotation}
``[Leech] dangled the problem under the nose of a few symmetry aficionados, Donald Coxeter among them, and he spread the word, telling his friend John McKay. McKay tempted Conway, and Conway eventually took the bait \ldots McKay was forever peddling mathematical unions, people with people, people with ideas.''
\end{quotation}

After Chilton, McKay's international life began.
In 1969, John crossed the Atlantic with wife Wendy, to be Visiting Assistant Professor at the California Institute of Technology.
There he continued to work on character tables, and formed a life-long friendship and collaboration with Hershy Kisilevsky, who would end up in Montreal with him some years later.
This brief sojourn in California was followed by McKay moving to Montreal in 1971 to take up a permanent post as Associate Professor in the School of Computer Science at McGill University.
The reason for moving to Canada, as John once said to me, was simple: he believed that the ratio between typical academic salary and cost-of-living was at the time maximised in Canada.
This could well be true.

McKay didn't get along with McGill, perhaps the place was too much establishment for him. His good friend and colleague, John Harnad, who visited him till the very last days, wrote in the obituary for the Canadian Centre de Recherches Math\'ematiques (CRM) \cite{HarnadObit}:
\begin{quotation}
``John McKay was a good friend, a distinguished colleague, and a very special person, with a unique character.
He was remarkable in many ways, both as an imaginative, creative mathematician, well-known for his incredible knack of discovering unexpected relations between seemingly unrelated fields, and as a clear-minded, perceptive thinker, with great honesty and integrity. It followed from this that he could not let things that he perceived as wrong to just go by, without probing their causes. He always did this with self-effacing good humour, in the spirit of inquiry rather than judgement, but sometimes with more persistence than certain academic administrators could much appreciate \ldots Besides his remarkable mathematical insights and contributions, it would be very incomplete to not say a bit more about John McKay’s human qualities, his very high ethical integrity, and his wonderfully wide scope of interests and knowledge, on historical, scientific, and cultural matters. In particular, he had a very deep sense of the rights and wrongs regarding how things are done in an academic institutional setting, and never ceased to question things that seemed to him as unjustified or inappropriately prioritised.
\end{quotation}

In 1974, John moved a few blocks away to Concordia University, where he remained until his final days, being promoted to full professor of computer science in 1979, and then to a joint full professor in mathematics and computer science in 1990, and Distinguished Emeritus Professor when he retired at the age of 70 in 2009.
A peregrination that paralleled his move to Concordia was when John and Wendy divorced in 1985, and while Wendy went on to marry fellow mathematician Robert Moody, John married Trinh Vo in 1988. Trinh had been by John's side with great affection and devotion.
She definitely organised his life. 
For a decade after John and I became friends and collaborators in 2010, we talked to each other almost every day on Skype \cite{mckaybio}.
This was greatly facilitated by the 5-hour time difference between Montreal and London/Oxford, because his mornings coincided well with my afternoons. Trinh would join us from time to time, especially when I had my two children, whom they saw growing up over the video chats. In many ways, McKay became grandpa John.

When it became too difficult for Trinh to take care of John after he had a fall -- I visited them in 2016 when I was giving a lecture at McGill and saw the precariousness of their stairs -- she arranged for John to be at Louis-Riel Nursing Home in 2017, where he remained until his last day. Skype became more difficult and Trinh tried whenever she could to get John onto an iPad so that he could continue to chat to me and to his other friends.
One rather poignant tradition Louis-Riel has is that a short obituary be written for every person at the time of admission.
Trinh asked me to do this and so more than five years ago I had the chance to rehearse the life story I'm writing now.
That brief homage to John was published by his alma mater in the Dulwich College Alleyne Club Yearbook \cite{dulwichObit}.
John died peacefully on 19 April 2022 with Trinh by his side.
Just a few days before, she managed to get me on the phone with him because, as she put it, my voice ``always cheered him up.''

McKay's life-long contributions to mathematics and to bridging ideas and people were duly celebrated and recognised when he was elected a Fellow of the Royal Society of Canada in 2000.
In 2003 he won the CRM-Fields Prize for Mathematics, Canada's highest honour for a mathematician.
In 2007 a conference was jointly organised by Concordia University and the Universit\'e de Montr\'eal to honour John and to marvel at what he brought to the world of mathematics. It was appropriately named
``Groups and Symmetries: From Neolithic Scots to John McKay'' \cite{mckayConf}.

%%%%%%%%%%%%%%%%%%%%%%%
\section{Mathematical Miracles}

\begin{quotation}
``John McKay was an exceptional mathematician. His broad mathematical interests and knowledge, together with his curiosity, intuition and insights, led to his discovery of important new connections between different areas of mathematics and the opening up of whole new areas of research.''\\
\flushright{- Leonard Soicher \cite{soicherObit}}
\end{quotation}

Einstein described the year 1905 as his {\it annus mirabilis} wherein he published four papers, each of which revolutionised our understanding of the universe.
In recounting McKay's contributions to mathematics it is also expedient to emphasise {\em his} three miracles, which transpired over the course of a decade, and in each case spurred  whole new research areas that remain active to this day.
They are (i) {\it Moonshine}, for which he was most famous, (ii) {\it McKay's A-D-E Correspondence}, which he considered his favourite, and (iii) the {McKay conjecture} on group representations, which was his first major observation.
Since John's most extraordinary ability was to link unimaginably disparate branches of mathematics, and since I deeply believe every graduate student should know and appreciate these miracles, I shall introduce them with sufficient generality and necessary context to at least attempt a conveyance of their beauty.

\begin{figure}[h]
\centerline{
(a) \includegraphics[trim=0mm 0mm 0mm 0mm, clip, width=2in]{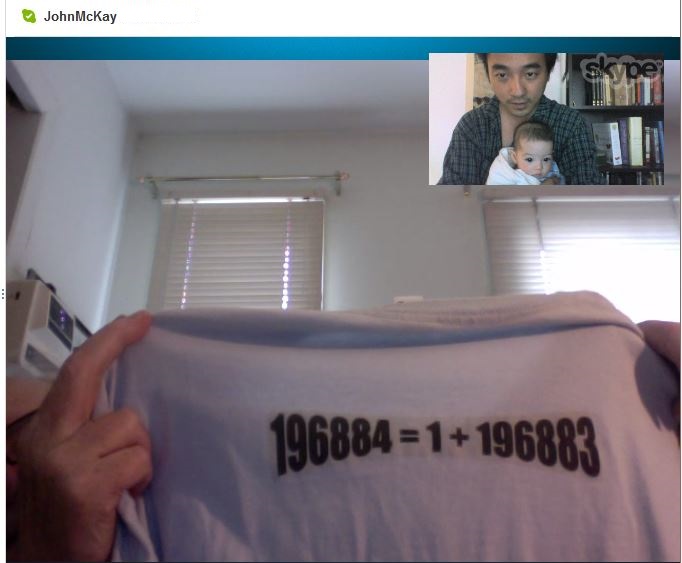}
(b) \includegraphics[trim=0mm 0mm 0mm 0mm, clip, width=2.2in]{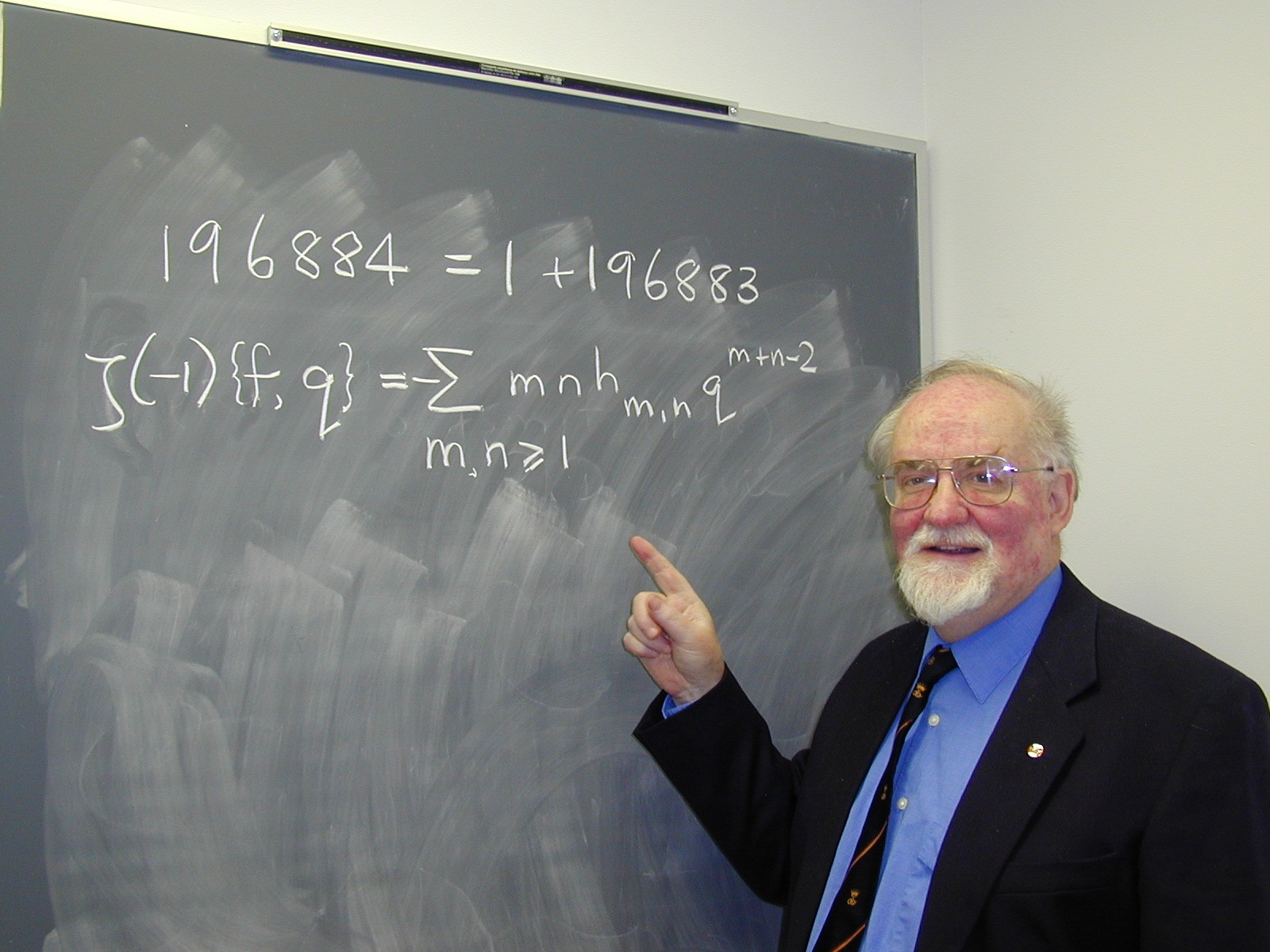}
}
\caption{(a) John showing the T-shirt to my then infant daughter Kitty, in a Skype call in 2014; 
(b) A well-known photo of McKay at the black-board, courtesy of John Harnad.} 
\label{f:196884}
\end{figure}

\subsection{``This is Moonshine!''}
McKay's most famous equation is
\begin{equation}\label{196884}
196884 = 196883 + 1 \ .
\end{equation}
The Fields Institute makes T-shirts with this printed on them and indeed next to John's ashes in {\it Cimeti\`ere Notre-Dame-des-Neiges} is the photo of John pointing to this equation (see Figure \ref{f:196884}).
What on earth does this equation mean?!
Let's explain the left and right in turn.

\def\sl{SL(2; \mathbb{Z})}

\subsubsection{The Modular Invariant}
Consider one of the most important concepts in mathematics, the modular group  $\sl := \left\{
\left( \begin{array}{cc}
a & b \\ c & d
\end{array} \right) : a,b,c,d \in \mathbb{Z}; ad - bc = 1
\right\}$, the group of $2 \times 2$ integer matrices with unit determinant\footnote{
Sometimes, the modular group is defined projectively as $PSL(2; \mathbb{Z})$ where the determinant $ad - bc = \pm 1$.
Much of our discussions carry through the $\pm$ sign.
}.
This group is ubiquitous, manifesting in combinatorics, geometry, number theory, and so on.
For one thing, the infinite number of parallelograms formed by the columns of its elements all have the same unit area as that of the identity element $\mathbb{I}_{2\times 2}$, which corresponds to the unit square.
For instance, the fact that the extremely skew one whose sides are generated by vectors $\{489, 100\}$ and $\{2000, 409\}$ has area 1 is highly non-trivial. 

Is there a function invariant under its action? To address this, we complexify the plane and recast the $\sl$ action as the integer linear fractional transformation $\tau \mapsto \frac{a \tau + b}{c \tau + d}$ for $\tau \in \mathbb{C}$, and seek an analytic function invariant under this transformation. Such a function would generalise well-known ones that are invariants under important subgroups of transformations, such as the trignometric function $\sin(2\pi \tau)$ being invariant under $\tau \mapsto \tau+1$.
The functions invariant {\em up to a weight factor} are the famous {\em modular forms}. The Eisenstein series prescribes them: for $k \in \mathbb{Z}_{>1}$, we have
\begin{equation}\label{eisenstein}
E_{2k}(\tau) = (2 \zeta(2k))^{-1} \sum\limits_{(0,0) \neq  (m,n) \in \mathbb{Z}^2} (m + n \tau)^{-2k} 
\quad \Rightarrow \quad
E_{2k}\left( \frac{a \tau + b}{c \tau + d} \right) = (c \tau+d)^{2k} E_{2k}(\tau) \ ,
\end{equation}
with the weight-factor $(c \tau+d)^{2k}$.
The series is absolutely convergent for $\tau \in \mathcal{H}$, the upper-half plane, where $\Im(\tau) > 0$.
The numerical factor $2 \zeta(2k)$ -- which involves the Riemann zeta-function evaluated at even positive integers, and thus can be expressed in terms of Bernoulli numbers -- has deep number-theoretical meaning which will not concern us presently.

Now, at the cost of introducing a pole, we can get rid of the weight factors, by defining
\begin{equation}
j(\tau) := 1728 \frac{E_4(\tau)^3}{E_4(\tau)^3 - 27 E_6(\tau)^2} 
\quad \Rightarrow \quad
j\left( \frac{a \tau + b}{c \tau + d} \right) = j(\tau)
\ .
\end{equation}
Incidentally, Ramanujan's familiarity with the $1728 = 12^3$ pre-factor of this function is surely the fundamental reason behind the Taxicab legend (the Lambert series of $E_4$ having the coefficient $n^3$).
Thus defined, it is clear that the properties of $E_{2k}$ under modular transformations from (\ref{eisenstein}) renders $j(\tau)$ invariant under $\sl$. However, it is {\em not} referred to as a ``modular form of weight zero'' because it is not holomorphic: unlike the modular forms $E_{2k}(\tau)$, $j(\tau)$ has a simple pole at $\tau = i \infty$.

This function $j$ is usually attributed to Klein \cite{klein} and is mostly called the Klein $j$-invariant.
However, McKay had always told me that it would have been known to at least Gauss and Eisenstein, so he insisted on it being simply called the {\em modular invariant function}, which emphasises the fact that it is {\em unique}, in the sense that it is the only function on the upper-half-plane invariant under $\sl$, which is holomorphic away from the simple pole at $i \infty$ and normalised so that
\begin{equation}
j(\exp(2 \pi i /3)) = 0 \ , \quad  j(i) = 1728 \ .
\end{equation}
Any other modular invariant (with different values at these points) is a rational function in $j$.

There are countless beautiful properties of $j(\tau)$, ranging from geometry (e.g., it being a detector of isomorphism classes of elliptic curves) to number theory (e.g., its connection to class field theory), but the one most pertinent to us is its series development. Defining $q := \exp(2 \pi i \tau)$, the Laurent series in $q$ (and equivalently, the Fourier series in $\tau$) is
\begin{equation}\label{seriesJ}
j(q) = q^{-1} + 744 + 196884 q + 21493760 q^2 + 864299970 q^3 + 20245856256 q^4 + \cdots
\end{equation}
We can see that the leading term of $q^{-1}$ signifies the simple pole at $q=0$, or $\tau = i \infty$.
All coefficients are integers, as can be seen from the integer Fourier expansion of the Eisenstein $E_{2k}$, but the meaning of these seemingly meaningless integers had to wait almost a century until McKay came along.

\subsubsection{Finite Simple Groups}
Now, for something completely different.
Parallel to analysis, geometry and number theory, another great development of the 19th century was the theory of groups, with its incipience as the tragic scribbles of Galois the night before his fateful duel.
Ever since Galois' idea of organising the permutations amongst the roots of polynomials -- in what we now call the Galois group -- the classification of finite groups had been a major theme in abstract algebra.

First, the problem boils down to the classification of finite {\em simple} groups, i.e., those which do not have a non-trivial proper normal subgroup.
We remind ourselves that a normal subgroup $N$ of a group $G$ is a subgroup closed under conjugation: $g n g^{-1} \in N$ for all $n\in N$ and all $g\in G$; 
this normality is denoted as $N \triangleleft G$. Normality permits one to form the quotient $G/N$.
Every finite group $G$ has a {\em composition series} which in a way encodes its structure: 
$\mathbb{I} = H_0 \triangleleft H_1 \triangleleft \ldots \triangleleft H_n = G$ where each $H_i$ is a maximal proper normal subgroup of $H_{i+1}$ and each quotient $H_{i+1}/H_i$ is simple.
The Jordan--H\"older theorem dictates that, up to permutations and isomorphisms, the composition series is unique.
The composition series for a finite group is like prime factorisation for integers, and thus simple groups to finite groups are like prime numbers to integers.

The classification of finite simple groups is itself a saga.
Galois noticed that the alternating group $\mathcal{A}_n$ on $n \geq 5$ elements was simple. This, together with what we now call Galois Theory of field extensions, was enough to resolve the long-standing problem of why generic polynomials of degree 5 and above have no solutions in radicals. In this sense, the classification problem began in 1832.
It is a remarkable journey spanning more than a century when the classification was officially completed.
A significant portion of the proof involved computer algebra, for which McKay was a pioneer.
The definitive volumes, consisting of almost ten thousand pages detailing the steps, are still an ongoing project \cite{finitesimple}.
For those uninitiated to pore over these tomes, there is a leisurely but technical account in \cite{gallian} which is very readable and even comes with a song at the end to help you remember the punchline:
\begin{theorem}
The finite simple groups are
\begin{enumerate}
\item Cyclic groups $\mathbb{Z}/p\mathbb{Z}$ of size $p$ prime;
\item The Alternating group $\mathcal{A}_n$ of even permutations on $n \geq 5$ elements;
\item The Lie groups defined over finite fields;
\item 26 Sporadic groups.
\end{enumerate}
\end{theorem}
Note that I am using ``size'' instead of ``group order'' and avoided using ``degree'' for the alternating groups. This is a nod to John, who always thought the words ``order'' and ``degree'' have been far too abused in mathematics.
In the above, Types (i-iii) are somewhat `classical' objects and come in infinite families.
Type (i) is clearly simple because they don't have any non-trivial proper subgroups.
Type (ii) is Galois' original observation.
Type (iii) is interesting in that Lie groups are usually associated with continuous symmetries, but defined over finite fields $\mathbb{F}_{p^r}$ they become discrete and finite; to Lie groups we will shortly return in the next section.
Type (iv), the sporadic groups\footnote{
There is the Tits group, which is sometimes called the 27th sporadic, and sometimes grouped together with the Lie-type family.
} are the most mysterious: they are all outliers and have rather large size -- the smallest of these, Mathieu's group $M_{11}$ discovered in 1861, is of size 7920.
The existence and construction of these 26 occupied a good part of the century-long story of the classification.
Aside from much active research they continue to engender, their {\em very being} remains puzzling. It is as if I told you that there are families of prime numbers like the Mersenne primes or Fermat primes, but then there are 26 specific primes of considerable size which are simply out there and are special and there are no more!

%%%
\subsubsection{Monstrous Moonshine}
The largest of the sporadics is of size
\begin{equation}
\begin{array}{c}
808017424794512875886459904961710757005754368000000000 = \\
2^{46} \cdot  3^{20}  \cdot 5^9 \cdot  7^6 \cdot  11^2 \cdot  13^3 \cdot  17 \cdot  19 \cdot  23 \cdot  29 \cdot  31 \cdot  41 \cdot  47 \cdot  59 \cdot  71 \ .
\end{array}
\end{equation}
The enormity of this number gave rise to the name Monster.
Its existence was first predicted independently by Bernd Fischer and Robert Griess \cite{Fischer-Griess} in the early 1970s.
In 1982, Griess published the first construction of the Monster, and a simplified construction was later published by
Conway.
John Thompson proved uniqueness, subject to the existence of a faithful
complex representation of dimension 196883, and Simon Norton showed that such a
representation must exist, subject to the structure of the prime-element
centralisers. A complete uniqueness proof for the Monster 
was finally published in 1989 by Griess, Meierfrankenfeld, and Segev \cite{GMS}.

This faithful complex representation of dimension $196883 = 47 \cdot 59 \cdot 71$ is fascinating. It is the smallest one.
In other words, the smallest matrices which generate the Monster -- in fact only two such matrices are needed, as all finite simple groups are generated by two matrices! -- are of size $196883 \times 196883$.
With today's technology, writing down these two matrices is not terrible, it requires about 1 Tb of disk-space. Writing down all the elements of the Monster is impossible -- it has more elements than the number of atoms of earth.
Conway summarises the discovery of the Monster well: ``There's never been any kind of explanation of why it's there, and it's obviously not there just by coincidence. It's got too many intriguing properties for it all to be just an accident.''

The legend goes as follows.
In 1978, John Conway\footnote{
Remember that Conway and McKay knew each other from the late 1960s when McKay dangled Leech's lattice to Conway.
It turns out that the Monster is intimately related to the automorphism group of the Leech lattice.
} showed the number 196883 to John McKay\footnote{At the time, this 196883-dimensional representation was still conjectural. Based on this conjectured degree, Fischer, Livingstone, and Thorne, in an amazing piece of work at the University of Birmingham, constructed the entire character table of the Monster. Thus the data were available for Thompson, Conway, and Norton to check that McKay's ideas were not so crazy!}
When McKay saw the number, his reply was Eq.~(\ref{196884}).
Conway replied that this is ridiculous, ``it's moonshine!'' Moonshine is English slang for (1) illegal home-brew and (2) something outrageously crazy. Since both Johns were English, they would get the reference.

Conway's outburst is well-justified.
The left hand-side of Eq.~(\ref{196884}) lives in the world of number theory and analytic functions; the right hand-side lives in the world of finite groups and representation theory.
The two worlds are so distant in the landscape of mathematics that it would indeed be crazy for them to be inter-related.
To make sure this is not just a numerical coincidence, they checked the next few (huge) coefficients
\begin{equation}\label{co-moon}
\begin{array}{rcl}
1 & = & r_1 \\
196884 & = & r_1 + r_2  \\
21493760 & = & r_1 + r_2 + r_3 \\
864299970 & = & 2r_1 + 2r_2 + r_3 + r_4 \\
20245856256 & = & 3r_1 + 3r_2 + r_3 + 2r_4 + r_5  = 2r_1+ 3r_2 + 2r_3 + r_4 + r_6\\
\ldots
\end{array}
\end{equation}
where on the left are the Fourier coefficients of $j(q)$ in Eq.~(\ref{seriesJ}) and $r_{n=1, \ldots, 194} = \{$ 1, 196883, 21296876, 842609326, 18538750076, 19360062527, \ldots, 258823477531055064045234375 \} are the dimensions of the irreducible representations of the Monster - there is a finite number, a total of 194, of them because it is a finite group.
Thus all the infinite number of integer coefficients of the modular invariant are encoded as simple linear combinations of the dimensions of the irreducible representations of the Monster sporadic simple group.
This can no longer be just a coincidence, especially because it involves such astronomical integers.

Thompson, who got the Fields Medal for his work on the simple groups and who discovered a sporadic that bears his name, suggested that there might exist a graded vector-space, now called the moonshine module, $V^\natural = \oplus_n V_n$ such that $j(q) - 744 = \sum_n {\rm Tr}(\mathbb{I} |_{V_n}) q^n$, because the trace of the identity $\mathbb{I}$ on a vector-space is its dimension\footnote{
The constant term 744 does not show up in Moonshine and needed to be subtracted away.
However, it is a remarkable fact that $744/3 = 248$ and has to do with the fact that $j(q)^{1/3}$ encodes the irreducible representations of the Lie group $E_8$ which we will encounter shortly \cite{kac}.
}.
This would then generalise to other elements $g$ - or conjugacy classes, since trace is conjugacy invariant - of the Monster as
\begin{equation}\label{M-T}
j_g(q) = \sum_n {\rm Tr}(g |_{V_n}) q^n \ . 
\end{equation}
This expansion is now called the McKay--Thompson series.
After Conway took the bait, he and Norton formalised his initial exclamation of incredulity \cite{ConwayNorton} into:
\begin{conjecture}  (The Moonshine Conjecture)
The McKay--Thompson series $j_g(q)$ is the unique modular invariant of a group lying between $\Gamma_0(N)$ and its normaliser $\Gamma_0(N)^+$ in the modular group $SL(2; \mathbb{Z})$.
\end{conjecture}
Here, $\Gamma_0(N)$ is the level $N$ congruence subgroup of the modular group consisting of the usual elements {\small $\left( \begin{array}{cc}
a & b \\ c & d
\end{array} \right)$}, $ad-bc=1$, but with the further congruence condition $c \bmod N \equiv 0$.
When $N=1$ and $g = \mathbb{I}$, $j_g(q)$ is the modular invariant $j(q)-744$ for the entire modular group, associated to McKay's original observation.
There is, of course, a lot more to the story.
The conjecture had many more precise details, and entire books have been devoted thereto, e.g., \cite{ronan,gannon,dusautoy} contain excellent accounts\footnote{For a short, enjoyable, but still technical, account, the reader is referred to \cite{Tatitscheff:2019okj}, by my former student Tatitscheff, who worked with John and me.}.
My intention here is merely to entice the reader with some highlights.
After immediate response from the community \cite{fong,FLM} in constructing the Monster module $V^\natural$, using ideas from conformal field theory and string theory, in particular vertex operator algebras, Borcherds proved the Moonshine Conjectures for the Monster, for which he was awarded the Fields Medal in 1998.

This surprising connection between modularity, finite groups and lattices, representation theory as well as theoretical physics, prompted by McKay's penetrating intuition, has become a field in and of itself, under the name of Moonshine.
Generalising beyond the Monster group to other sporadics, lattices beyond the Leech lattice, and $j(q)$ to partition functions, remains an active field of research: since the 1990s, there have been some 300 papers on this subject on the ArXiv.

When Conway turned 80, there was a conference and volume dedicated to him \cite{conwayVol}. McKay and I wrote the opening chapter \cite{MR3682584}, a one-page mathematical (half) joke entitled ``Moonshine and the Meaning of Life'' where we summed the squares of the first 24 coefficients of $j(q)$, as well as for the related Ramanujan tau-function, and reduced modulo 70  (both 24 and 70 being key to the Leech lattice) and obtained 42. Conway was very amused.
I am honoured that the very last published work I have with John was my editing and tidying up a set of notes of his final in-person lecture on his thoughts and new ideas about Moonshine \cite{MR4484467}.

Let me conclude this section with  the words of Norton \cite{nortonbio}:
\begin{quotation}
``I can explain what Monstrous Moonshine is in one sentence, it is the voice of God.''
\end{quotation}

%%%%%%%%%%
\subsection{A-D-E}
Our second mathematical miracle is what John always told me to be his favourite, even more than Moonshine.
Again, this is the bridging of unrelated fields: finite groups, graphs, geometric singularities and Lie theory. It is now called the {\em McKay Correspondence} and is widely believed to be but the tip of the iceberg for many more deep connections in mathematics yet to be discovered.

\subsubsection{Lie Groups}
Sophus Lie's original ambitious programme was to classify continuous symmetries by considering the group of transformations on the variables of a multi-variate analytic function. This led to the concept which we now call {\em Lie groups}, i.e., continuous groups which are also manifolds. With many standard examples we are familiar.
For example, the circle $S^1$ is a group by identifying its elements as unimodular complex numbers $\exp(i \theta)$, under  usual multiplication; this is the Lie group $U(1)$.
Another example is $SU(2)$, the group of $2 \times 2$ complex matrices $M$ such that $M^H M = \mathbb{I}_{2\times2}$ and $\det(M) = 1$, where $H$ is the conjugate transpose.
In other words, $M$ are the special unitary $2 \times2$ matrices, forming the Lie group which geometrically is the 3-sphere $S^3$.
Yet another example, of which we will shortly make use, is $SO(3)$, the group of $3 \times 3$ real matrices $M$ such that $M^TM = \mathbb{I}_{3\times3}$ and $\det(M) = 1$, i.e., special orthogonal matrices; this Lie group is topologically the real projective space $\mathbb{R}\mathbb{P}^3$.
The above-mentioned are examples of {\em classical Lie groups}, which are more or less familiar to the 19th century mathematician\footnote{
To reconnect with the previous section, these continuous Lie groups, when defined not over the reals or the complex numbers, but over finite fields, become discrete, and furnish the Lie family of finite simple groups.}.

Lie's insight was that the Lie group, being a manifold, has a tangent space at the group identity. This tangent space is a vector space endowed with a binary operation inherited from the group law that renders it an algebra; the algebra is now dubbed the {\em Lie algebra} and the binary operation, the Lie bracket.
This group-algebra correspondence transforms the study of a Lie group to that of the Lie algebra whose vector space structure renders it more manageable.

The brilliant works of Killing, Cartan, Levi et al.~subsequently reduced the classification of Lie groups to that of Lie algebras as complex vector spaces. In analogy to the finite groups, Lie algebras can be built up from the {\em simple Lie algebras}, which are the ones that have no nonzero proper ideals. The Levi structure theorem dictates that any Lie algebra, when quotiented by its radical (largest solvable ideal), is a direct sum of simple ones. Thus again, simple Lie algebras are the ``primes'' and the complete list is, rather prosaically dubbed, as follows:
\begin{theorem}\label{thm:Lie}
The simple Lie algebras over $\mathbb{C}$ are
\begin{enumerate}
\item Type $A_n$: $\mathfrak{sl}_{n+1}(\mathbb{C})$;
\item Type $B_n$: $\mathfrak{so}_{2n+1}(\mathbb{C})$;
\item Type $C_n$: $\mathfrak{sp}_{2n}(\mathbb{C})$;
\item Type $D_n$: $\mathfrak{so}_{2n}(\mathbb{C})$;
\item The Exceptionals: $E_6$, $E_7$, $E_8$, $F_4$, and $G_2$.
\end{enumerate}
\end{theorem}
In the above $n \in \mathbb{Z}_{>1}$ and Types A to D are called the Classical Lie algebras which come in infinite families, while Type E to G, are called the Exceptionals. One could see an analogy to the classification of simple finite groups!

Now, the correspondence between Lie algebras and Lie groups, for explicit matrices, is the exponentiation map.
Thus, Type A are traceless $(n+1) \times (n+1)$ complex matrices, which, when exponentiated, give rise to $(n+1) \times (n+1)$ unit-determinant matrices, i.e., $SL(n; \mathbb{C})$. 
The real form is $\mathfrak{su}_{n+1}$ -- in the sense that $\mathfrak{sl}_{n+1}(\mathbb{C}) = \mathfrak{su}_{n+1} \otimes_{\mathbb{R}} \mathbb{C}$ --
and are traceless, anti-Hermitian matrices.
These exponentiate to special unitary matrices.
In other words, Type A Lie algebra corresponds to the Lie group $SU(n+1)$.
Likewise, Type B and D Lie groups are the special orthogonal groups of odd and even dimension respectively;
Type C is the symplectic group.
The Exceptionals are more mysterious, with the root systems living in dimensions 6, 7, 8, 4, and 2.

%%% https://math.stackexchange.com/questions/3275149/how-are-sun-sln-and-mathfraksln-mathbbc-related

%%%%%%%%%%%%%%%
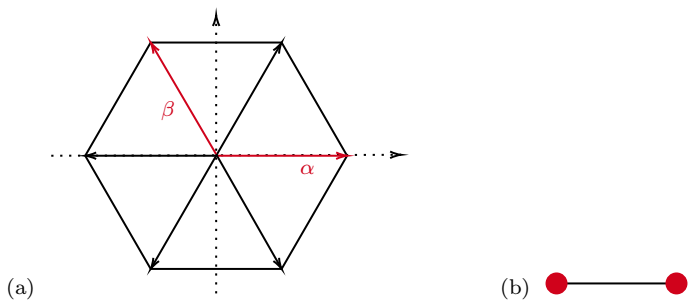
\begin{figure}
\centerline{
(a) \tikzset{every picture/.style={line width=0.75pt}} %set default line width to 0.75pt        

\begin{tikzpicture}[x=0.75pt,y=0.75pt,yscale=-1,xscale=1,scale=0.50]

%Shape: Regular Polygon [id:dp7588037763474962] 
\draw   (352,149) -- (286.5,262.45) -- (155.5,262.45) -- (90,149) -- (155.5,35.55) -- (286.5,35.55) -- cycle ;
%Straight Lines [id:da2310805804424818] 
\draw [fill={rgb, 255:red, 74; green, 144; blue, 226 }  ,fill opacity=1 ] [dash pattern={on 0.84pt off 2.51pt}]  (56,149) -- (406,148.07) ;
\draw [shift={(406,148.07)}, rotate = 180] [color={rgb, 255:red, 0; green, 0; blue, 0 }  ][line width=0.75]    (10.93,-3.29) .. controls (6.95,-1.4) and (3.31,-0.3) .. (0,0) .. controls (3.31,0.3) and (6.95,1.4) .. (10.93,3.29)   ;
%Straight Lines [id:da04271542693755759] 
\draw  [dash pattern={on 0.84pt off 2.51pt}]  (221,287) -- (221,10) ;
\draw [shift={(221,8)}, rotate = 90] [color={rgb, 255:red, 0; green, 0; blue, 0 }  ][line width=0.75]    (10.93,-3.29) .. controls (6.95,-1.4) and (3.31,-0.3) .. (0,0) .. controls (3.31,0.3) and (6.95,1.4) .. (10.93,3.29)   ;
%Straight Lines [id:da37358938950122766] 
\draw [color={rgb, 255:red, 208; green, 2; blue, 27 }  ,draw opacity=1 ][fill={rgb, 255:red, 208; green, 2; blue, 27 }  ,fill opacity=1 ]   (221,149) -- (330,149) -- (350,149) ;
\draw [shift={(352,149)}, rotate = 180] [color={rgb, 255:red, 208; green, 2; blue, 27 }  ,draw opacity=1 ][line width=0.75]    (10.93,-3.29) .. controls (6.95,-1.4) and (3.31,-0.3) .. (0,0) .. controls (3.31,0.3) and (6.95,1.4) .. (10.93,3.29)   ;
%Straight Lines [id:da1894021771177017] 
\draw    (221,149) -- (285.5,37.28) ;
\draw [shift={(286.5,35.55)}, rotate = 120] [color={rgb, 255:red, 0; green, 0; blue, 0 }  ][line width=0.75]    (10.93,-3.29) .. controls (6.95,-1.4) and (3.31,-0.3) .. (0,0) .. controls (3.31,0.3) and (6.95,1.4) .. (10.93,3.29)   ;
%Straight Lines [id:da8217098681466009] 
\draw [color={rgb, 255:red, 208; green, 2; blue, 27 }  ,draw opacity=1 ]   (221,147.5) -- (156.51,37.28) ;
\draw [shift={(155.5,35.55)}, rotate = 59.67] [color={rgb, 255:red, 208; green, 2; blue, 27 }  ,draw opacity=1 ][line width=0.75]    (10.93,-3.29) .. controls (6.95,-1.4) and (3.31,-0.3) .. (0,0) .. controls (3.31,0.3) and (6.95,1.4) .. (10.93,3.29)   ;
%Straight Lines [id:da12942909895907584] 
\draw    (221,149) -- (92,149) ;
\draw [shift={(90,149)}, rotate = 360] [color={rgb, 255:red, 0; green, 0; blue, 0 }  ][line width=0.75]    (10.93,-3.29) .. controls (6.95,-1.4) and (3.31,-0.3) .. (0,0) .. controls (3.31,0.3) and (6.95,1.4) .. (10.93,3.29)   ;
%Straight Lines [id:da17150732947665293] 
\draw    (221,149) -- (156.5,260.72) ;
\draw [shift={(155.5,262.45)}, rotate = 300] [color={rgb, 255:red, 0; green, 0; blue, 0 }  ][line width=0.75]    (10.93,-3.29) .. controls (6.95,-1.4) and (3.31,-0.3) .. (0,0) .. controls (3.31,0.3) and (6.95,1.4) .. (10.93,3.29)   ;
%Straight Lines [id:da04049632378761636] 
\draw    (221,147.5) -- (285.51,260.71) ;
\draw [shift={(286.5,262.45)}, rotate = 240.32] [color={rgb, 255:red, 0; green, 0; blue, 0 }  ][line width=0.75]    (10.93,-3.29) .. controls (6.95,-1.4) and (3.31,-0.3) .. (0,0) .. controls (3.31,0.3) and (6.95,1.4) .. (10.93,3.29)   ;

% Text Node
\draw (302,155.4) node [anchor=north west][inner sep=0.75pt]  [color={rgb, 255:red, 208; green, 2; blue, 27 }  ,opacity=1 ]  {$\alpha $};
% Text Node
\draw (163,91.4) node [anchor=north west][inner sep=0.75pt]  [color={rgb, 255:red, 208; green, 2; blue, 27 }  ,opacity=1 ]  {$\beta $};

\end{tikzpicture}\qquad\qquad
(b)
\tikzset{every picture/.style={line width=0.75pt}} %set default line width to 0.75pt
\begin{tikzpicture}[x=0.75pt,y=0.75pt,yscale=-1,xscale=1]
%Straight Lines [id:da20852789495313062] 
\draw    (135,125) -- (195,125) ;
%Shape: Circle [id:dp39452525982288655] 
\draw  [color={rgb, 255:red, 208; green, 2; blue, 27 }  ,draw opacity=1 ][fill={rgb, 255:red, 208; green, 2; blue, 27 }  ,fill opacity=1 ] (130,125) .. controls (130,122.24) and (132.24,120) .. (135,120) .. controls (137.76,120) and (140,122.24) .. (140,125) .. controls (140,127.76) and (137.76,130) .. (135,130) .. controls (132.24,130) and (130,127.76) .. (130,125) -- cycle ;
%Shape: Circle [id:dp3761818363378251] 
\draw  [color={rgb, 255:red, 208; green, 2; blue, 27 }  ,draw opacity=1 ][fill={rgb, 255:red, 208; green, 2; blue, 27 }  ,fill opacity=1 ] (190,125) .. controls (190,122.24) and (192.24,120) .. (195,120) .. controls (197.76,120) and (200,122.24) .. (200,125) .. controls (200,127.76) and (197.76,130) .. (195,130) .. controls (192.24,130) and (190,127.76) .. (190,125) -- cycle ;
\end{tikzpicture}
}
\caption{(a)The root system and (b) Dynkin diagram of $A_2$.}
\label{f:root}
\end{figure}
%%%%%%%%%%%%%%%%%%

Finally, Dynkin and Coxeter gave a marvellous graphical representation of Theorem \ref{thm:Lie}, by further reducing the problem of studying simple algebras to {\em root systems}, using the so-called Cartan basis which we will not delve into now.
Suffice to say that the classification of Lie groups reduces first to that of simple Lie groups, thence to that of simple Lie algebras, and finally to that of positive simple roots in a root system:
\begin{definition}
A root system $\Phi$ is a finite collection of non-zero vectors, called roots, in Euclidean space $\mathbb{R}^r$ with the standard dot product $(\cdot, \cdot)$ such that
\begin{itemize}
\item The roots span $\mathbb{R}^r$;  $r$ is called the rank;
\item The only scalar multiples of a root $\alpha \in \Phi$ that are still in $\Phi$ are $\pm \alpha$ and none other;
\item For any two roots $\alpha, \beta$, the vector $\beta - 2 \frac{(\alpha, \beta)}{(\alpha, \alpha)} \alpha$, which is the reflection through the hyperplane perpendicular to $\alpha$, is also a root \footnote{
	For Lie algebras, the factor $2 \frac{(\alpha, \beta)}{(\alpha, \alpha)}$ in the reflection is an integer.
}
\end{itemize}
In $\Phi$, we can always choose, potentially in many ways, the set $\Phi^+$ of positive roots which are such that
\begin{itemize}
\item For each root $\alpha \in \Phi$, only one of $\pm \alpha$ is in $\Phi^+$;
\item For any two distinct roots $\alpha, \beta \in \Phi^+$ such that $\alpha + \beta$ is a root, then $\alpha + \beta \in \Phi^+$.
\end{itemize}
Finally, a simple root is a positive root which cannot be written as the sum of two other positive simple roots.
We denote the set of simple positive roots as $\Delta \subset \Phi^+ \subset \Phi$; the number of simple roots turns out to be equal to the rank $r$.
\end{definition}
The archetypical example of a root system in $\mathbb{R}^2$ is given in Fig.~\ref{f:root} (a), prescribed by the vertices of a regular hexagon. One can check that the roots are $\{\pm \alpha, \pm \beta, \pm(\alpha+\beta)\}$, the positive roots are $\{ \alpha, \beta, \alpha+\beta \}$ and the positive simple roots are $\{\alpha, \beta\}$, which subtend an angle of $\frac23 \pi$. This root system is called $A_2$, and it does indeed correspond to the simple Lie algebra $A_2$, and the associated Lie group $SU(3)$.

The root system representation can be even further simplified, using a Coxeter--Dynkin diagram (or just Dynkin diagram), using the following rules:
(1) Let $\Delta$ be the simple roots and form a finite graph with each node corresponding to a simple root -- therefore there will be $r$ nodes where $r$ is the rank.
(2) The edges between two nodes are determined by the angle between the two corresponding simple roots using the rule:
\[
\mbox{
$\pi / 2$: no edge; 
$2 \pi / 3$: single edge; 
$3 \pi / 4$: double edge; 
$5 \pi / 6$: triple edge.
}
\]
It turns out that these are the only possible angles.
(3) We introduce direction, i.e., turn the edges into arrows - by letting the longer simple root point to the shorter one; if the roots are of the same length, then there's no need for direction.
In the above example of $A_2$, there are therefore two nodes corresponding to $\alpha$ and $\beta$, there is a single edge between them because they are at 120 degrees, and they are also of equal length. The Dynkin diagram is in Fig.~\ref{f:root} (b).

\begin{figure}
\centerline{
\includegraphics[trim=0mm 0mm 0mm 0mm, clip, width=4in]{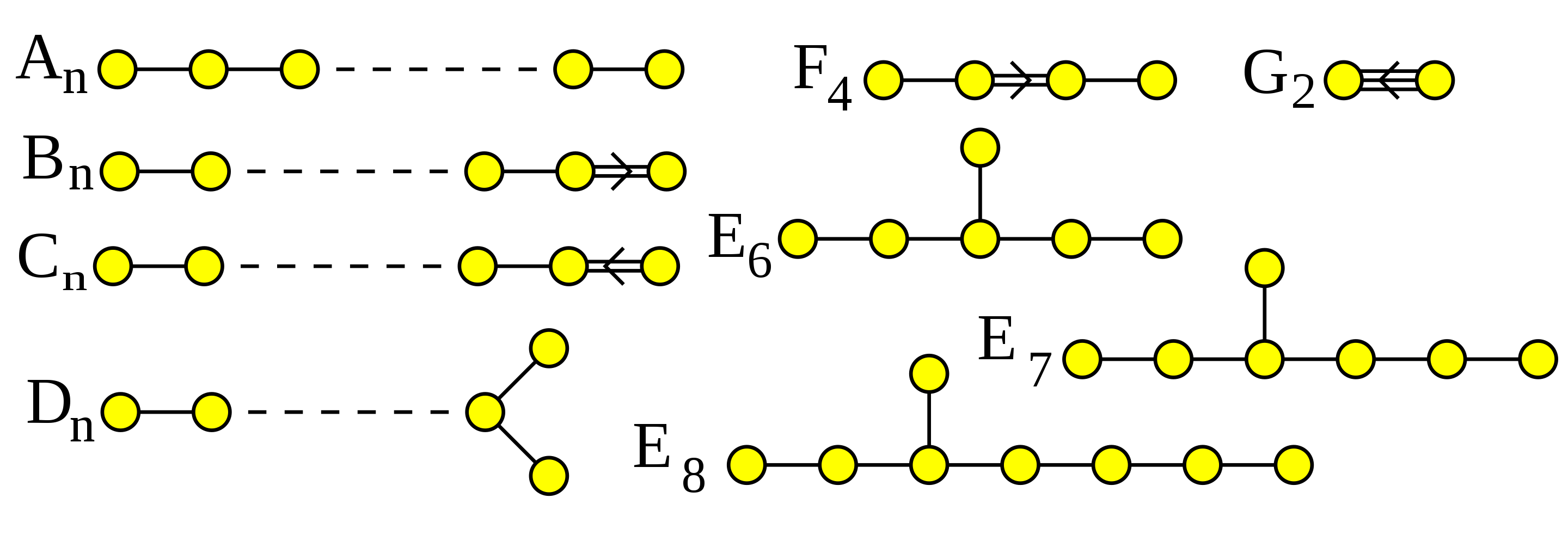}
}
\caption{The Dynkin diagrams of the simple Lie algebras/groups (source: Wiki).}
\label{f:Dynkin}
\end{figure}

In summary, the graphical representation of Theorem \ref{thm:Lie} is Fig.~\ref{f:Dynkin}.
One readily notices that of these diagrams, Types A-D-E have only single (and undirected) edges, meaning that all the simple roots are of equal length, and are either perpendicular, or at an angle of $2\pi/3$ to each other.
These are called {\em simply-laced} Dynkin diagrams, or simply {\bf A-D-E diagrams}.
They come in two infinite families $A_n$ and $D_n$, as well as three exceptionals $E_{6,7,8}$.

%%
%%Graphics3D[{Opacity[0.8],   PolyhedronData["Icosahedron", "GraphicsComplex"]}, Boxed -> False]
\begin{figure}[h!]
\[
\begin{array}{|c|c|c|}\hline
\mbox{Tetrahedron} &
\mbox{Cube \qquad Octahedron}&
\mbox{Dodecahedron\qquad Icosahedron}
\\
\includegraphics[trim=0mm 0mm 0mm 0mm, height=1in]{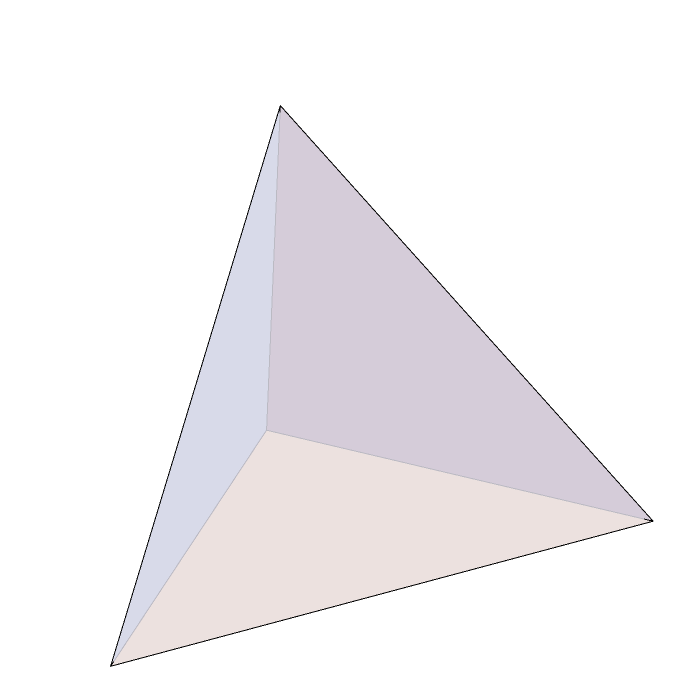}
&
\includegraphics[trim=0mm 0mm 0mm 0mm, height=1in]{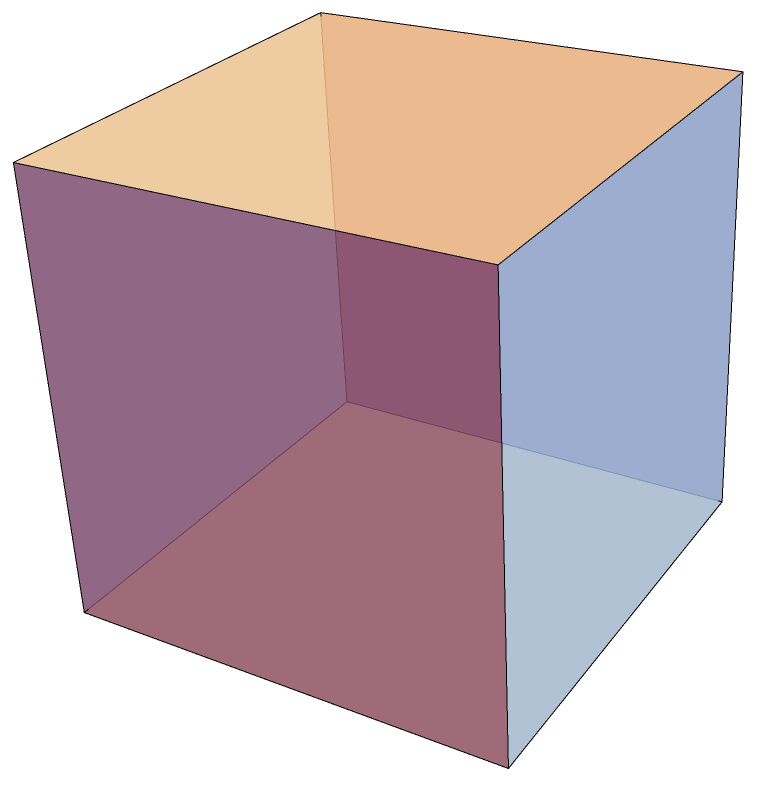}
\includegraphics[trim=0mm 0mm 0mm 0mm, height=1in]{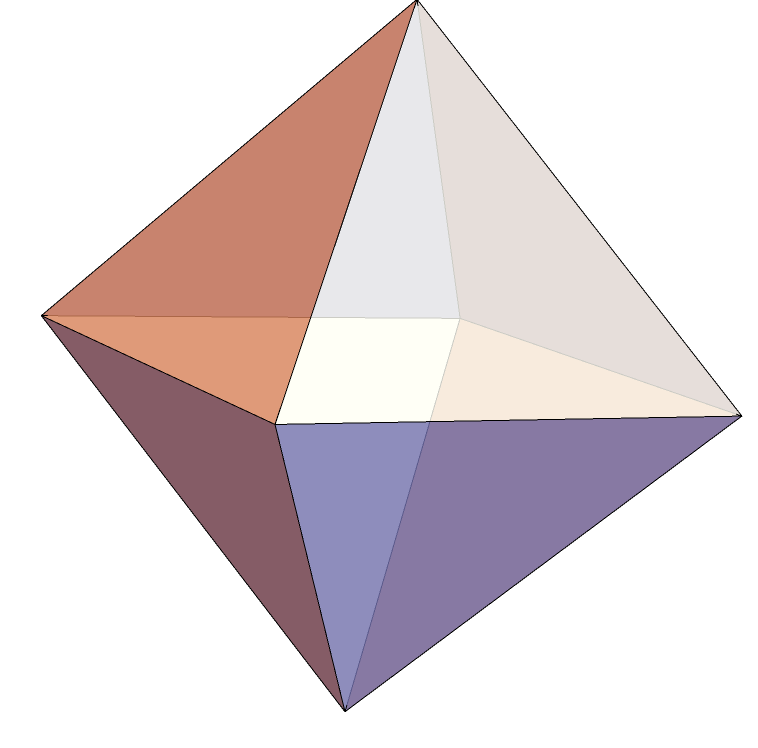}
&
\includegraphics[trim=0mm 0mm 0mm 0mm, height=1in]{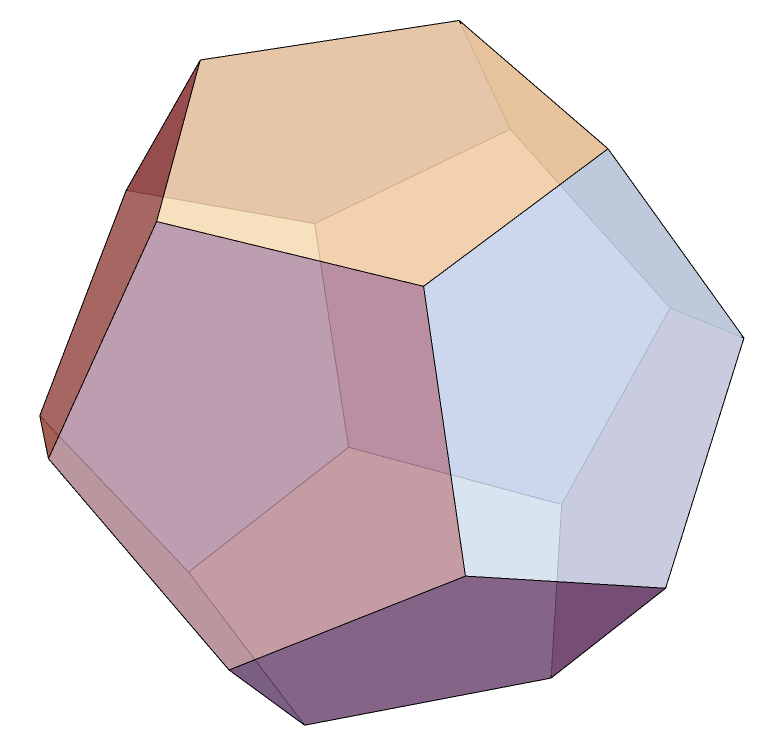}
\includegraphics[trim=0mm 0mm 0mm 0mm, height=1in]{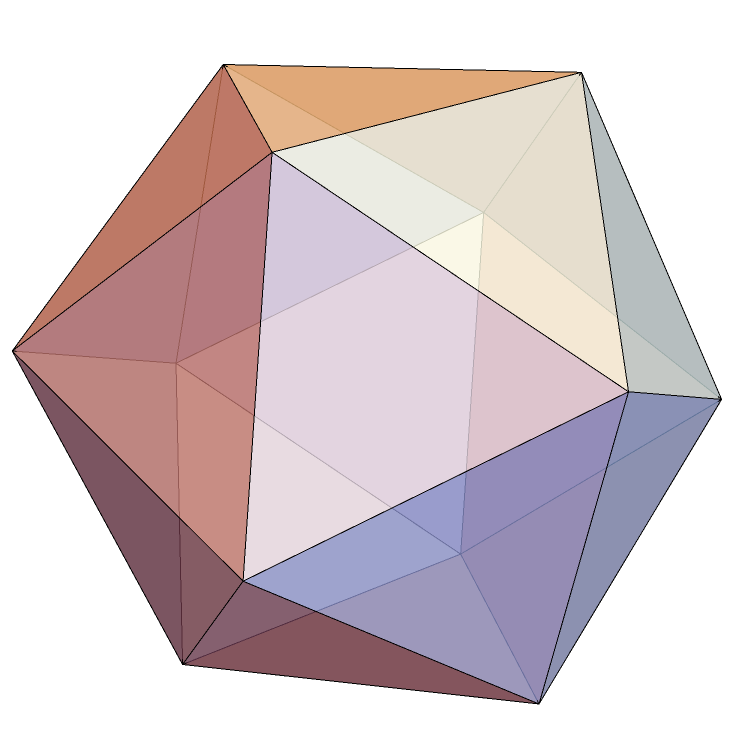}
\\
\hline
\end{array}
\]
\caption{The five Platonic solids: the regular tetrahedron, cube, octahedron, dodecahedron and icosahedron.
We have grouped them according to duality.}
\label{f:Plato}
\end{figure}

\subsubsection{Platonic Solids}
Now, for something completely different.
It is an old question what are the discrete symmetries of Euclidean space.
Phrased in a way ancient Greeks would understand: what are the regular convex shapes?
They are clearly (1) the regular polygons and (2) the five Platonic solids, viz., the regular tetrahedron, the cube, the regular octahedron, the regular dodecahedron and the regular icosahedron. For reference, we show these five famous shapes in Fig.~\ref{f:Plato}.

Phrased in a more modern way: what are the discrete finite subgroups of $SO(3)$, the group of special orthogonal matrices that constitute the rotational group of symmetries of $\mathbb{R}^3$?
These are the finite groups that are the symmetries of the aforementioned regular convex shapes (polytopes).

The regular $n$-gons are degenerate in the sense that they lie in the plane, rather than being properly 3-dimensional. For them
the symmetry groups are (1) the cyclic group of size $n$ generated by $2\pi / n$ rotations; and (2) the dihedral group of size $2n$, generated by $2\pi / n$ rotations and the flip along any axis of symmetry.
For the five platonic solids there are only three symmetry groups. This is because of (graph) duality, in the sense of replacing faces with vertices and edges with perpendicular edges.
The cube--octahedron are dual to each other, and likewise for the dodecahedron--icosahedron, while the tetrahedron is self-dual.
Hence, there are three symmetry groups:
\begin{itemize}
\item (3a) the tetrahedral group which turns out to be the alternating group $\mathcal{A}_4$, of size $4!/2 = 12$, 
\item (3b) the cube--octahedral group (or simply the octahedral group) which is the symmetric group $\mathcal{S}_4$, of size $4!=24$, and 
\item (3c) the dodeca--icosahedral group (or simply the icosahedral group) which is the alternating group $\mathcal{A}_5$, of size $5!/2 = 60$.
\end{itemize}
Everything I've said here would have been known to the 19th century mathematician, and knowledge of the Lie groups and algebras, to those of the early 20th.
Yet, once again, it took waiting for many decades until John came along to notice an extraordinary link between this set of symmetry groups of $\mathbb{R}^3$: 2 infinite families and 3 outliers, and the seemingly utterly unrelated set of simply-laced Dynkin diagrams: 2 infinite families and 3 exceptionals.
In 1979, McKay announced at a conference on groups, his observation which came to be known as the {\em McKay Correspondence}.

\subsubsection{McKay Correspondence}
John said to me several times that his favourite and proudest observation -- in his gentle humility he never used words such as `discovery' when referring to himself -- was this 1979 correspondence.
McKay first introduces a complex representation by lifting $SO(3)$ to its double cover $SU(2)$ (incidentally, these two Lie groups share the same Lie algebra). This lift is familiar to physicists because it is a key difference between classical (vector) and quantum-mechanical (spinor) angular momentum. 
On a group level, the discrete symmetry groups are called {\em binary} and are double covers of what was mentioned above. Specifically, we have (note that the cyclic groups, being Abelian, have no non-trivial double cover)
\begin{theorem}\label{su2}
The finite discrete subgroups of $SU(2)$ are:
\[
\begin{array}{|c|c|c|} \hline
\mbox{Name} & \mbox{Presentation} & |G| \\ \hline \hline
\mbox{Cyclic}
	& \gen{R \ | \ R^n = \mathbb{I}} \simeq \mathbb{Z}/n\mathbb{Z} & n  \\ \hline
\mbox{Binary Dihedral}
    & \gen{R,S \ | \ R^n = S^2 = (RS)^2}  & 4n  \\ \hline
\mbox{Binary Tetrahedral}
	& \gen{R,S,T \ | \ RST=R^2 = S^3 = T^3} & 24  \\ \hline
\mbox{Binary Octahedral}
	& \gen{R,S,T \ | \ RST=R^2 = S^3 = T^4} & 48  \\ \hline
\mbox{Binary Icosahedral}
	& \gen{R,S,T \ | \ RST=R^2 = S^3 = T^5 } & 120  \\
\hline
\end{array}
\]
\end{theorem}
Next, McKay takes the defining complex 2-dimensional representation ${\bf R_2}$ (which necessarily exists because our groups live in $SU(2)$) of each of these groups, and forms the tensor decomposition over all the irreps ${\bf r_i}$ 
\begin{equation}\label{mckay2}
{\bf R_2} \otimes {\bf r_i} = \bigoplus\limits_j a_{ij} {\bf r_j}
\end{equation}
where $a_{ij} \in \mathbb{Z}_{\geq 0}$ denote the multiplicity.
While such tensor decompositions for groups, discrete or continuous, are standard, somehow it occurred to no one that this should be done to the $SU(2)$ subgroups.

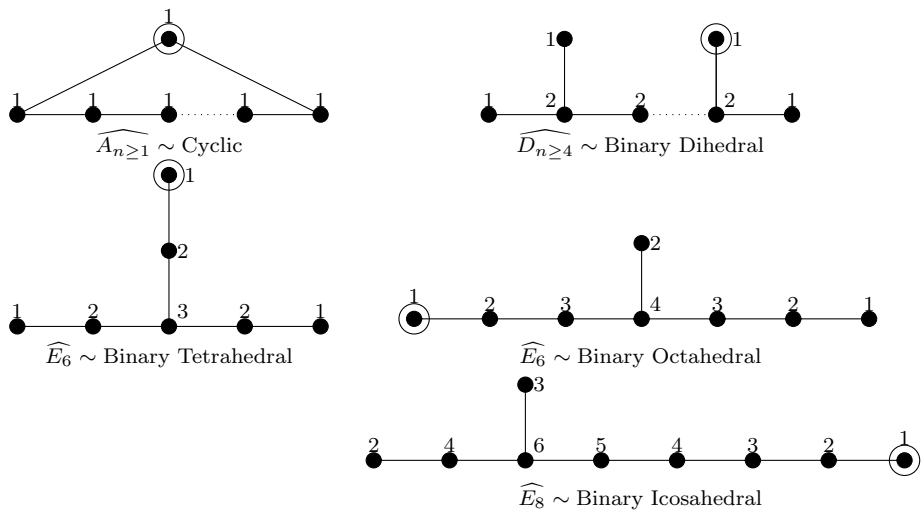
\begin{figure}
\[
\begin{array}{cc}
\begin{tikzpicture}
\draw[fill=black]
(0,0)
	circle [radius=.1] node [above] {1} --
(1,0)
	circle [radius=.1] node [above] {1} --
(2,0)
	circle [radius=.1] node [above] {1}
(3,0)
	circle [radius=.1] node [above] {1} --
(4,0)
	circle [radius=.1] node [above] {1}
;
\draw[dotted]
(2,0) -- (3,0)
;
\draw[fill=black]
(2,1)
	circle [radius=.1]
;
\draw
(2,1)
	circle [radius=.2] node [above=1mm] {1}
;
\draw[solid]
(0,0) -- (2,1)
;
\draw[solid]
(4,0) -- (2,1)
;
\end{tikzpicture}

&%%%%%%%

\begin{tikzpicture}
\draw[fill=black]
(0,0)
	circle [radius=.1] node [above] {1} --
(1,0)
	circle [radius=.1] node [above left] {2} --
(2,0)
	circle [radius=.1] node [above] {2}
(3,0)
	circle [radius=.1] node [above right] {2} --
(4,0)
	circle [radius=.1] node [above] {1}
(1,1)
	circle [radius=.1] node [left] {1}
(3,1)
	circle [radius=.1] node [left] {}
;
\draw[solid]
(3,0) -- (3,1)
;
\draw[dotted]
(2,0) -- (3,0)
;
\draw[solid]
(1,0) -- (1,1)
;
\draw[solid]
(3,0) -- (3,1)
;
\draw
(3,1)
	circle [radius=.2] node [right=1mm] {1}
;
\end{tikzpicture}

\\

\widehat{A_{n\geq 1}} \sim \mbox{Cyclic}
&
\widehat{D_{n \geq 4}} \sim \mbox{Binary Dihedral}

\\

\begin{tikzpicture}
\draw[fill=black]
(0,0)
	circle [radius=.1] node [above] {1} --
(1,0)
	circle [radius=.1] node [above] {2} --
(2,0)
	circle [radius=.1] node [above right] {3} --
(3,0)
	circle [radius=.1] node [above] {2} --
(4,0)
	circle [radius=.1] node [above] {1}
(2,1)
	circle [radius=.1] node [right] {2}
(2,2)
	circle [radius=.1] node [] {}
;
\draw[solid]
(2,0) -- (2,1)
;
\draw[solid]
(2,1) -- (2,2)
;
\draw
(2,2)
	circle [radius=.2] node [right=1mm] {1}
;
\end{tikzpicture}

&

\begin{tikzpicture}
\draw[fill=black]
(-1,0)
	circle [radius=.1] node [above] {} --
(0,0)
	circle [radius=.1] node [above] {2} --
(1,0)
	circle [radius=.1] node [above] {3} --
(2,0)
	circle [radius=.1] node [above right] {4} --
(3,0)
	circle [radius=.1] node [above] {3} --
(4,0)
	circle [radius=.1] node [above] {2} --
(5,0)
	circle [radius=.1] node [above] {1}
(2,1)
	circle [radius=.1] node [right] {2}
;
\draw[solid]
(2,0) -- (2,1)
;
\draw
(-1,0)
	circle [radius=.2] node [above=1mm] {1}
;
\end{tikzpicture}

\\

\widehat{E_6} \sim \mbox{Binary Tetrahedral}
&
\widehat{E_6} \sim \mbox{Binary Octahedral}

\\

&
\begin{tikzpicture}
\draw[fill=black]
(0,0)
	circle [radius=.1] node [above] {2} --
(1,0)
	circle [radius=.1] node [above] {4} --
(2,0)
	circle [radius=.1] node [above right] {6} --
(3,0)
	circle [radius=.1] node [above] {5} --
(4,0)
	circle [radius=.1] node [above] {4} --
(5,0)
	circle [radius=.1] node [above] {3} --
(6,0)
	circle [radius=.1] node [above] {2} --
(7,0)
	circle [radius=.1] node [] {}
(2,1)
	circle [radius=.1] node [right] {3}
;
\draw[solid]
(2,0) -- (2,1)
;
\draw
(7,0)
	circle [radius=.2] node [above=1mm] {1}
;
\end{tikzpicture}
\\
&
\widehat{E_8} \sim \mbox{Binary Icosahedral}

\end{array}
\]
\caption{The McKay Correspondence between the affine simply-laced Dynkin diagrams and the discrete finite subgroups of $SU(2)$, or between the affine Lie algebras and the Platonic solids.}
\label{f:affineADE}
\end{figure}

Finally, McKay recognises the matrices $a_{ij}$, which happen to be symmetric and with 2 on the diagonal and 0 or 1 otherwise!
Explicitly, they can be expressed \cite{Hanany:1998sd} in terms of the character table $\chi_{\gamma}^{(i)}$ of the finite group $G$ as 
\begin{equation}\label{aij}
a_{ij} = \frac{1}{|G|} \sum\limits_{\gamma=1}^r C_\gamma \chi_\gamma^{{\bf R_2}} \chi_\gamma^{(i)} \overline{\chi}_\gamma^{(j)}
\ , 
\end{equation}
where the sum is over the $r$ conjugacy classes and $C_\gamma$ is the size of the $\gamma$-th conjugacy class.
He considers them to be the adjacency matrices of some finite graphs -- now called the McKay Quivers -- by associating a node to each row/column and an edge between vertices $i$ and $j$ whenever $a_{ij}=1$.
The result is presented in Fig.~\ref{f:affineADE}.
These are almost the simply-laced ADE Dynkin diagrams, just with an extra node!
It turns out that these are called the {\em affine} ADE Dynkin diagrams and encode an infinite generalisation of the corresponding Lie algebras called affine Lie algebras.
The extra affine node, in the McKay Correspondence, is associated to the identity element in the finite group.

Once again, John managed to link two fields of mathematics which had not talked to each other for at least half a century, Lie theory and symmetry groups of the Platonic solids.
Perhaps mathematicians had missed this because both sides needed some manipulation: affinisation on the Lie side and lift to $SU(2)$ on the rotational groups side; in addition, one needed the genius\footnote{To demystify all this a little, what inspired John was, first, the pattern of 2 infinite families and 3 outliers, and second, that the sum of squares of certain Coxeter labels in the Dynkin diagrams (which I will not discuss here for sake of space) gives rise to the size of the finite groups.} of Eq.~(\ref{mckay2}).
One might wonder whether the Moonshine sporadics and the McKay Correspondence exceptionals might be related.
The initial seeds were indeed planted by John and elaborated in \cite{He:2015yoa}, awaiting further clarification.

My very first paper and a good part of my PhD thesis was to find generalisations of the McKay Correspondence and especially to see their realisations in superstring theory in physics. Thus from the very outset of my career, I worshipped McKay.
To complete the circle, the very last thing on which I worked with him was a textbook \cite{ADEbook}, together with Cameron and Dechant, aimed at a wide mathematical audience, entitled ``A-D-E''. Tragically, John will never see the completion of this project.

The McKay Correspondence did much more than bridge two fields. It prompted mathematicians to look for the ADE meta-pattern everywhere. It seems ubiquitous across mathematical disciplines ranging from algebraic geometry to conformal field theory, and has gained almost a religious status. V.~Arnol'd found some 20 disparate instances of what he appropriately called the ``Trinity'', $E_{6,7,8}$ \cite{Arnold}.
If you ever run into a classification problem in mathematics where there are 2 infinite families and 3 exceptional cases, then beware! Chances are, they fall into an A-D-E type of meta-pattern.

%
% #### Binary Icosahedral
%	f := FreeGroup( "r", "s", "t");
%	g := f / [(f.1*f.2*f.3)^(-1)*f.1^2, (f.1*f.2*f.3)^(-1)*f.2^3, (f.1*f.2*f.3)^(-1)*f.3^5];
%	Size(g);
%
% #### p = 2
%	s := SylowSubgroup(g, 2); n := Normalizer(g, s);

%%%%%%%%%%%%
\subsection{Characters \& Primes}

We finally proceed to John's earliest miracle \cite{McKayCon,MR286904}, made between 1971-2, which still remains one of the most important open problems in representation theory (for a recent review see \cite{evseev}).
This is commonly known as the McKay Conjecture, later generalised and formalised \cite{alperin} to the Alperin-McKay Conjecture.
In the words of Markus Linckelmann, 
\begin{quotation}
``The McKay conjecture, from 1972, and its block theoretic generalisation known
as the Alperin-McKay conjecture, from 1976, are instances of a host of
on-the-surface miraculous numerical coincidences that haunt modular
representation theory of finite groups.''
\end{quotation}

Once again, it begins with a miraculous observation by John, made shortly after he completed his PhD.
Consider the complex characters $\chi_\gamma^{(i)}$ of a finite group $G$ -- the computer calculation of which was the subject of his thesis.
Recall that characters are traces of the matrix representations and  $\chi_\gamma^{(i)}$ is an $r \times r$ square table where $r$ is equal to both the number of conjugacy classes of $G$ and the number of irreducible representations, termed `irreps'. The rows, indexed by $\gamma$, go through the irreps, and the first row is always just all 1s, corresponding to the trivial representation where all group elements are represented by 1. That is, $\chi_1^{(i)} = 1$.
The columns, indexed by $i$, go through the conjugacy classes, and the first column is always the list of dimensions of the irreps, being the trace of the identity matrix representation of $\mathbb{I}_G$. That is, $\chi_\gamma^{(1)}$ are the dimensions of the irreps, like the $r_i$ in Eq.~(\ref{co-moon}) for the Monster.
We made use, in the foregoing, of the character table for both the McKay-Thompson series in Eq.~(\ref{M-T}) and in inverting the McKay decomposition in Eq.~(\ref{mckay2}).

Now, fix a prime $p$ and consider a Sylow $p$-subgroup $S_p$ of $G$. 
Recall that a $p$-subgroup is a subgroup of size which is a prime-power $p^n$, and a Sylow $p$-subgroup is maximal in the sense that it is not a proper subgroup of any other $p$-subgroup.
Next, consider the normaliser of $S_p$ in $G$, i.e., $N_G(S_p) := \{ g \in G : gS_pg^{-1} = S_p \}$.
The McKay Conjecture, which remains open for arbitrary finite groups\footnote{
Interestingly, John's original ``McKay Conjecture'' was for finite simple groups and the prime $p=2$, which is now proved \cite{Ruh}).
}, states
\begin{conjecture} (McKay)
Denote by $Irr_{p'}(H)$ the set of characters such that $p$ {\em does not} divide the dimension of the irreps of a group $H$, i.e., $p \nmid \chi_\gamma^{(1)}$. Then
\[
| Irr_{p'}(G) | = | Irr_{p'}(N_G(S_p)) | \ .
\]
That is, the two sets have the same cardinality.
\label{conj:mckay}
\end{conjecture}

Let us see this conjecture in action.
It should work for any finite group, so we might as well take one with which we are familiar.
From Theorem \ref{su2}, consider the binary icosahedral group (which turns out to be isomorphic to $SL(2;5)$, the 2-dimensional special linear group defined over the finite field $\mathbb{F}_5$).
The character table of this group is presented in Table \ref{t:char} (a).
In the fine tradition of McKay, one could readily use computer algebra \cite{GAP4} to obtain this.
We see that there are $r=9$ irreps and conjugacy classes, and we have added an extra row on top to denote the sizes of the conjugacy classes for completeness.
We can check that this row sums to 120, the size of the group and, by standard Schur orthgonality for character tables, that the square sum of the first column $1^2 + 2^2 + 2^2 + 3^2 + 3^2 + 4^2 + 4^2 + 5^2 + 6^2 = 120$ also\footnote{Indeed, on an unrelated note, one can further check that substituting this character table into Eq.(\ref{aij}) gives the adjacency matrix of the affine $E_8$ Dynkin diagram. As one can check that the character table in part (b) gives the affine $E_6$ Dynkin diagram.}.

%(FromGapForm[ReadList["ct.temp", String]][[1]] /.    w[5] -> Exp[2 Pi I/5] //    FullSimplify) /. {1/2 (1 + Sqrt[5]) -> \[Phi],    1/2 (1 - Sqrt[5]) -> phibar, 1/2 (-1 + Sqrt[5]) -> -phibar,    1/2 (-1 - Sqrt[5]) -> -\[Phi]} // TeXForm
%%
\begin{table}[h!]
(a)
$
\begin{array}{|c|c|c|c|c|c|c|c|c|c|} \hline
1 & 12 & 12 & 20 & 30 & 20 & 12 & 1 & 12 \\ \hline \hline
 1 & 1 & 1 & 1 & 1 & 1 & 1 & 1 & 1 \\ \hline
 2 & -\phi'& -\phi  & 1 & 0 & -1 & \phi'& -2 & \phi  \\ \hline
 2 & -\phi  & -\phi'& 1 & 0 & -1 & \phi  & -2 & \phi'\\ \hline
 3 & \phi  & \phi'& 0 & -1 & 0 & \phi  & 3 & \phi'\\ \hline
 3 & \phi'& \phi  & 0 & -1 & 0 & \phi'& 3 & \phi  \\ \hline
 4 & -1 & -1 & 1 & 0 & 1 & -1 & 4 & -1 \\ \hline
 4 & -1 & -1 & -1 & 0 & 1 & 1 & -4 & 1 \\ \hline
 5 & 0 & 0 & -1 & 1 & -1 & 0 & 5 & 0 \\ \hline
 6 & 1 & 1 & 0 & 0 & 0 & -1 & -6 & -1 \\ \hline
 \end{array}
$
(b)
$
\begin{array}{|c|c|c|c|c|c|c|}\hline
 1 & 4 & 4 & 6 & 1 & 4 & 4 \\  \hline \hline
 1 & 1 & 1 & 1 & 1 & 1 & 1 \\ \hline
 1 & \omega_3^2 & \omega_3^2 & 1 & 1 & \omega_3 & \omega_3 \\ \hline
 1 & \omega_3 & \omega_3 & 1 & 1 & \omega_3^2 & \omega_3^2 \\ \hline
 2 & -1 & 1 & 0 & -2 & 1 & -1 \\ \hline
 2 & -\omega_3 & \omega_3 & 0 & -2 & \omega_3^2 & -\omega_3^2 \\ \hline
 2 & -\omega_3^2 & \omega_3^2 & 0 & -2 & \omega_3 & -\omega_3 \\ \hline
 3 & 0 & 0 & -1 & 3 & 0 & 0 \\ \hline
\end{array}
$
\caption{The character tables of (a) the binary icosahedral group and (b) the normalizer of its Sylow 2-group.
Here, $\phi = \frac{1+\sqrt{5}}{2}$, $\phi' = \frac{1-\sqrt{5}}{2}$ and $\omega_3 = \exp(2 \pi i / 3)$.
}
\label{t:char}
\end{table}

Now let's pick the prime $p=2$.
The Sylow 2-group $S_2$ turns out to be the quaternion group $Q_8$ of size 8.
Its normaliser $N_G(Q_8)$ turns out to be $SL(2;3)$, which happens to be the binary tetrahedral group of size 24, whose character table is given in Table \ref{t:char} (b).
Here, there are $r=7$ irreps and conjugacy classes, and again we add a top row for the sizes of the conjugacy classes.

The characters for $G$ for which the entries in the first column are not divisible by $p=2$ are the rows which begin with $\{1,3,3,5\}$ from Table \ref{t:char} (a).
Likewise, from Table \ref{t:char} (b), the characters for $N_G(Q_8)$ for which the entries of the first column are not divisible by $2$ are the rows beginning with $\{1,1,1,3\}$.
Indeed, the cardinality of both sets is 4, as the conjecture demands.

There are collaborative efforts attempting to prove the conjectures using the classification of finite simple groups.
For $p=2$, Conjecture \ref{conj:mckay} was proved by Malle and Sp\"ath \cite{MS} in 2016, and the Alperin-McKay generalisation, again for $p=2$, was only proved last year by Ruhstorfer \cite{Ruh}, a week before John died.
Continuing to quote Linckelmann,
\begin{quotation}
``In parallel to the attempts using the classification of finite simple groups,
these numerical conjectures have generated a lot of work searching for structural
explanations in terms of module categories, their derived categories, drawing
on fields further afar, such as homotopy theory.
Apart for some special cases, structural explanations remain largely elusive --
translating  mystery into insight seems still out of reach.''
\end{quotation}

%%%%%%%%%%%%%%%%%%%%%%
\section{Legacy}

\begin{quotation}
``Aside from his energy and broad curiosity, the thing that struck me most was his phenomenal instinct for finding the singular mathematical pearl in a field of marbles. I think that what characterised his remarkable contributions to disparate mathematical areas, was his incredible intuition for identifying those (seemingly mundane) facts which were keys to unlocking deep underlying structures. These drew extensively on his ability to adapt computational methods to abstract mathematical issues \ldots
I will always be grateful for having had John as a friend and colleague. He made Concordia and the world a more interesting place.''
\flushright{ - Hershy Kisilevsky}
\end{quotation}

Officially, McKay, who was never into formalities, had only one PhD student, K.-C.~Young from McGill in 1974.
Yet he was enormously generous with his time and shared his ideas gregariously, especially with the younger generation.
Therefore, many of us feel as if we are in his intellectual family.
These have, over the years, included:
Leonard Soicher, a Master's student of John's, who on his recommendation went on to pursue his PhD with Conway, and then returned to be McKay's postdoc;
Abdellah Sebbar who was his postdoc;
Abdelkrim Elbasraoui who was a student at Concordia;
Roland Friedrich who worked with him on free probability theory and who called him regularly.
Additionally, there is myself, my students James Read and Valdo Tatitscheff, and many more.
We all felt like John's students.
Over his career spanning half a century, McKay had many influential publications, which I have collated here \cite{MR217029} - \cite{ADEbook}.
Of course, of all these, his three miracles distinguish themselves, and to their account I have devoted this obituary.

I have heard some non-mathematicians, and even an occasional mathematician, dismiss what I have described here.
They presume that if you make a conjecture a day, surely you will hit on something important sooner or later.
To these pompous and tiresome people, I mumble, ``Well, why don't you try it!'' 
It has been a tremendous privilege to have done what they could not, and that is to have worked with John McKay on somewhat of an equal footing.
To have had a glimpse of McKay's methods and thoughts has made me realise that to hit on one of his miracles is already astounding, but to do so thrice requires super-human powers.
I should also point out that everything John did was decades before the internet, when one could not simply google at will.
Today, if you typed `196884' into OEIS \cite{oeis}, you would find $j$-invariant coefficients and moonshine-related material.

John had nothing but his vast knowledge and perspicacious insight.
His mind was divine.

~\\
~\\

\begin{acknowledgements}\label{ackref}
I am grateful to John's wife Trinh McKay-Vo, his sister Liz Hurt, and his ex-wife Wendy McKay for all their tirelessly answering my questions on his life. Above all, Trinh, together with John, had very much become a part of my family; to her this obituary is dedicated. Sincere thanks also go to John's friends Leonard Soicher, John Harnad, Hershy Kisilevsky, and John Ball for providing various snippets of his life over different periods, as well as Madeleine Hall, John Harnad, Leonard Soicher, Pierre-Philippe Dechant, Thomas Hodgkinson, and Markus Linckelmann for many helpful comments.
\end{acknowledgements}

%
%
%
%
% Important: Do not put any empty line here.%
\affiliationthree{% in this example, two authors share an institution
  Yang-Hui He\\
  London Institute for Mathematical Sciences,
  Royal Institution of Great Britain,
  21 Albemarle Street, London W1S 4BS, UK \&
  Merton College, Oxford University, OX14JD, UK
   \email{hey@maths.ox.ac.uk}}
% Important: Do not put any empty line here.%
\end{document}